# Computational Study on the Impact of Gasoline-Ethanol Blending on Autoignition and Soot/NO$_x$ Emissions under Gasoline Compression Ignition Conditions

**Krishna C. Kalvakala**
Department of Mechanical and Industrial Engineering, University of Illinois at Chicago, IL, USA

**Harsimran Singh**
Department of Mechanical and Industrial Engineering, University of Illinois at Chicago, IL, USA

**Pinaki Pal***
Argonne National Laboratory, IL, USA

**Jorge Pulpeiro Gonzalez**
Argonne National Laboratory, IL, USA

**Christopher P. Kolodziej**
Argonne National Laboratory, IL, USA

**Suresh K. Aggarwal**
Department of Mechanical and Industrial Engineering, University of Illinois at Chicago, IL, USA

*Corresponding Author: pal@anl.gov*

***Abstract:*** *In the present work, computational fluid dynamics (CFD) simulations of a single-cylinder gasoline compression ignition (GCI) engine are performed to investigate the impact of gasoline-ethanol blending on autoignition, nitrogen oxide (NO$_x$), and soot emissions under low-load conditions. In order to represent the test gasoline (RD5-87), a four-component toluene primary reference fuel (TPRF) + ethanol (ETPRF) surrogate (with 10% ethanol by volume; E10) is employed. A three-dimensional (3D) engine CFD model employing finite-rate chemistry with a skeletal kinetic mechanism (including NO$_x$ sub-mechanism), adaptive mesh refinement (AMR), and hybrid method of moments (HMOM) is adopted to capture the in-cylinder combustion phenomena and soot/NO$_x$ emissions. The engine CFD model is validated against experimental data for three gasoline-ethanol blends: E10, E30 and E100, with varying ethanol content by volume. Model validation is carried out for a broad range of start-of-injection (SOI) timings (-21, -27, -36, and -45 crank angle degrees (CAD) after top-dead-center (aTDC)) with respect to in-cylinder pressure, heat release rate, combustion phasing, NO$_x$ and soot emissions. For relatively later injection timings (-21 and -27°aTDC), E30 yields higher amount of soot than E10; while the trend reverses for early injection cases (-36 and -45°aTDC). On the other hand, E100 yields the lowest amount of soot among all fuels irrespective of SOI timing. Further, E10 shows a non-monotonic trend in soot emissions with SOI timing: SOI-36>SOI-45>SOI-21>SOI-27, while soot emissions from E30 exhibit monotonic decrease with advancing SOI timing. NOx emissions from various fuels follow a trend of E10>E30>E100. On the other hand, NO$_x$ emissions increase as SOI timing is advanced for all fuels, with an anomaly for E10 and E100 where NO$_x$ decreases when SOI is advanced beyond -36°aTDC. Detailed analysis of the numerical results is performed to investigate the soot/NO$_x$ emission trends and elucidate the impact of chemical composition and physical properties on autoignition and emissions characteristics.*

## 1. INTRODUCTION

The growing energy demand in the transportation sector makes it necessary to build efficient and low-emission energy systems. Internal Combustion Engines (ICEs) are the most popular energy systems within the transportation sector and are predicted to contribute for meeting nearly 70% of transportation needs even in 2050 [1]. There has been an enormous amount of effort involved over last few decades in developing advanced combustion engines (ACEs) and in identifying suitable alternative/renewable fuels. These, when implemented together can render the energy generation process (for transportation) significantly efficient i.e., to maximize performance and minimize emissions. However, in order to achieve this overarching goal, there is a need to investigate the relationships between fuel composition (and associated properties) and its effect on combustion performance and emissions.

Gasoline Compression Ignition (GCI) mode is an advanced compression ignition (ACI) technology operating under low-temperature regime, which has gained significant interest in recent times. In this mode, gasoline is directly injected into the combustion chamber and ignited through compression. Several studies have shown the potential of GCI mode for achieving higher efficiencies and lower emissions than conventional diesel compression ignition (CI) engine [2-5]. Moreover, comprehensive investigations on the effect of various parameters, such as intake charge temperature, start of fuel injection, fuel injection pressure, and temperature, number of fuel injections, compression ratio, piston bowl geometry, and fuel reactivity were carried out to improve the performance of GCI mode of operation [5-12]. While the majority of these studies were focused on medium to high load conditions, there are limited studies which investigated the performance of GCI engines under low load conditions. The GCI engine under low load suffers from combustion instability which is mainly associated with poor reactivity of gasoline under low temperatures [13-17]. Weall et al. [17] investigated the impact of increasing intake pressure and temperature on the performance of low load GCI. It was observed that higher temperatures and higher pressures allowed for advanced crank angle corresponding to 50% total heat release (CA50). However, the higher temperatures resulted in substantial increase in $NO_x$ emissions. Roberts et al. [18] made similar observations about using higher intake temperatures and pressures. In addition, it was also observed that addition of a small amount of diesel to gasoline improved combustion stability at low load conditions. However, increasing fuel reactivity with addition of diesel significantly increased soot and $NO_x$ emissions. Vallinayagam et al. [19] computationally investigated various low-octane gasolines with varying levels of octane sensitivity ($S$) and distillation curves. The results indicated that the physical properties of the fuel have small impact on heat release rate characteristics under different levels of fuel stratification. Further, the use of gasoline with high $S$ was shown to improve the stability of GCI engine under low load conditions. However, the observations about $S$ are case specific considering the small sample of fuels considered. Zhang et al. [16] computationally examined the effect of fuel composition, chemical properties - Research Octane Number (RON), $S$, and physical properties on combustion phasing and emission characteristics. It was observed that fuels with lower RON improved reactivity while on the flip side produced higher soot emissions. Further, it was also observed that density has a significant effect on combustion phasing, soot, and $NO_x$ emissions while other physical properties had a very small impact on the performance. Higher density deteriorated air entrainment into the spray resulting in slower fuel-air mixing process, thereby producing higher soot. A recent study by Kalvakala et al. [20] showed that fuels with similar reactivities (in terms of CA10 and CA50) can exhibit significantly different sooting tendencies. The sooting tendency was strongly coupled to both fuel chemistry and physical properties with increased susceptibility to physical properties (mainly heat of vaporization (HoV) and viscosity) at advanced start of injection (SOI) conditions. On the other hand, the ignition delay time of these fuels was primarily driven by chemistry and showed very small dependency on physical properties. Further, the sooting tendency depicted a non-monotonic behavior with advancing SOI: -36>-45>-21>-27°aTDC. This is in contrast to fundamental understanding that higher mixing times available at advanced SOI timings would result in leaner fuel mixtures, thus lower soot emissions. The higher soot emissions at advanced SOI timings were a consequence of fuel's physical properties which increased the susceptibility of wall film formation (fuel deposits on piston and liner walls). These wall films resulted in pool fires leading to higher soot emissions. This indicates the paramount importance of fuel reactivity and properties on the performance of GCI mode of operation.

In the context of improving the reactivity and lowering emissions from GCI engines, blends of gasoline and biofuels were examined by several studies [21-33]. Ethanol has been identified as one of the potential low-carbon biofuels due to its higher volatility that allows faster mixing, thus improving the homogeneity of the cylinder charge [34]. Further, the higher latent heat of vaporization of ethanol can aid in lowering combustion temperatures, which results in suppressing the formation of $NO_x$. Further, it is hypothesized that the oxygen content in ethanol aids in suppressing soot emissions. Several fundamental and detailed engine studies have shown the potential benefits associated of blending gasoline with ethanol [20,28,35-38]. Presently, RD5-87 gasoline containing 10% (by volume) ethanol, with anti-knocking index (AKI) of 87, is adopted as the baseline fuel in the market of several economies across the globe. Woo et al. [36] compared the combustion and emission characteristics of diesel fuel and pure ethanol under CI engine conditions. Ethanol showed nearly 50% higher fuel conversion efficiency, despite having 36% lower calorific value than diesel. Although aforementioned studies attempted to examine the effect of blending ethanol with gasoline, the blending was limited to 20% by volume. Moreover, there are limited studies which systematically investigated the effect of blending ethanol in different proportions with gasoline on combustion and $NO_x$/soot emission characteristics of GCI engines operating under low load conditions.

In this context, the present study has three novel contributions. First, a robust computational fluid dynamic (CFD) model is developed that can predict the in-cylinder characteristics of GCI engine fueled by three different blends of gasoline/ethanol. The computational model is validated against experimental data generated as part of this study. Second, the effects of SOI timing on autoignition propensity and soot emissions from gasoline-like fuels under GCI low load conditions are investigated using validated CFD model. Lastly, the significance of ethanol concentration on autoignition propensity and soot emissions is systematically examined. For the same, three fuels – E10, E30, and E100 containing 10%, 30%, and 100% (by volume) ethanol in gasoline are considered. Moreover, the computational model is leveraged to isolate the impact of physical properties and chemical reactivity of fuels on autoignition propensity, $NO_x$ and soot emissions. The rest of the paper is organized as follows. Sections 2, 3 and 4 describe the experimental setup, operating conditions, and details about the engine CFD model. Subsequently, the Results and Discussion section is divided into four subsections. The first part presents the validation of chemical kinetic mechanism and CFD model against experimental data. The second part discusses the impact of fuel composition on autoignition and combustion phasing. The third and final part includes the discussion about the impact of fuel composition on soot and $NO_x$ emissions, respectively. The paper concludes with a summary of key findings.

## 2. EXPERIMENTAL SETUP

Figure 1 shows the Heavy-Duty Caterpillar (CAT) 3401 single-cylinder oil test engine (SCOTE) used for experiments. The engine is located at Argonne National Laboratory, where experiments for three different fuels of interest were performed. The engine has the standard 3401 SCOTE piston and is equipped with common rail direct injection system along with Next Gen injector. Table 1 provides the specifications of the engine. Soot mass was measured using an AVL Micro Soot Sensor. Additional details of the engine system are provided in a recent publication by the authors [20].

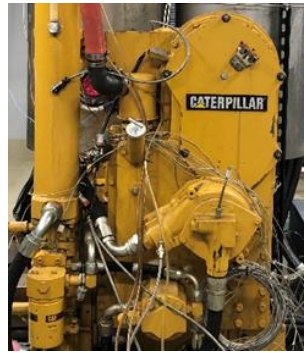

Figure 1: Experimental setup

Table 1: Specifications of the engine and fuel injection system.

| **Base Engine** | **Caterpillar SCE** |
|---|---|
| Injection System | DI Common Rail |
| Bore x Stroke | 137 × 165 mm |
| Displacement | 2.44 liters |
| Compression Ratio | 15.3:1 |
| Inj. Nozzle | 7x185 μm, 130° spray inclusion angle |

## 3. FUELS AND OPERATING CONDITIONS

Table 2 shows the composition of fuels investigated in the present study. Baseline gasoline used in the experiments is 87AKI E10 RD-587 gasoline. The second blend E30 is prepared by additional ethanol into RD5-87 such that the total ethanol in the fuel is 30% (by volume). The third fuel of interest is E100 containing 100% ethanol. Considering the complex chemical composition of gasoline used in experiments, a four-component surrogate was formulated which represents RON, *S,* and ethanol content in RD5-87 gasoline, and is adopted for the computational study. Also, the liquid density and heat of vaporization (HoV) of the surrogate are nearly similar to that of RD5-87The blend is generated using the artificial neural network (ANN) model developed at Lawrence Livermore National Laboratory (LLNL) [39,40]. The surrogate for the E30 blend is generated by adding ethanol to E10 surrogate such that the total amount of ethanol is 30% by volume.

Table 2: Composition and chemical properties of fuels considered.

| Fuel | RD5-87[1] Gasoline (vol%) | E10TPRF (E10) (vol%) | E30TPRF (E30) (vol%) | E100TPRF (E100) (vol%) |
|---|---|---|---|---|
| Isoparafinns/$IC_8H_{18}$ | 65.1 | 39.08 | 30.39 | 0.0 |
| Paraffins/$NC_7H_{16}$ |  | 21.03 | 16.35 | 0.0 |
| Aromatics/$C_6H_5CH_3$ | 20.1 | 29.91 | 23.26 | 0.0 |
| Ethanol | 9.9 | 9.98 | 30.0 | 100 |
| RON | 91.8 | 92.3 | 101.5[41] | 107.2[42] |
| *S* | 7.9 | 7.6 | 12.4[41] | 13.4[42] |
| Density (g/cc) | 0.754 | 0.748 | 0.763 | 0.785 |
| HoV (kJ/kg) | 380 | 416 | 593 | 874 |

The engine is operated at a global equivalence ratio (φ) of 0.32 and at a speed of 1200 rpm. The intake pressure and temperature at inlet valve opening are 1 bar and 145°C, respectively. The liquid fuel is injected using a 7-hole injector at an injection pressure of 500 bar. The compression ratio (CR) used in experiments was 15.26±0.45. The CR used in simulations is estimated by matching the motored pressure traces from simulations and experiments and fixed CR at 15.6:1 for all fuels and SOI timings. Four SOI timings are investigated: -21, -27, -36, and -45°aTDC in both experiments and CFD simulations. The injection timing used in CFD simulations is adjusted by 5 crank angle degrees (CAD) such that it accounts for the delay associated with the injection system adopted in the present experiments. This level of uncertainty is consistent with other recent experimental studies [7, 43, 44]. The experimental data presented in the current study is an average of 300 cycles. However, the CFD simulations performed are closed cycle simulations. In other words, they start at crank angles corresponding to inlet valve closing (IVC) and end at exhaust valve opening (EVO). The pressure at IVC (-169°aTDC) for simulations is the average pressure (of 300 cycles) at corresponding crank angle from experiments. On the other hand, the temperature at IVC is estimated using the Redlich-Kwong equation of state [45]. Moreover, in order to account for any residual gases from previous cycles, the CFD simulations consider the charge at IVC to be a mixture of air (95%) and residuals ($CO$, $CO_2$, $H_2O$, and $NO$) (5% by mass). The composition of residual gases used in simulations is an average over 300 experimental engine cycles. It must be noted that the mass of fuel injected for E10 fuel is ~47.5 mg while for other fuels is estimated such that the global phi ~~air-to-fuel ratio~~ is fixed at 0.32. In other words, the total fuel energy injected was held constant for all fuels under consideration. Further, the duration of injection (DOI) is adjusted such that the injection pressure is fixed at 500 bar.

---

[1] In addition to saturates, aromatics and ethanol, 4.9% (by vol) of the RD-87 gasoline content was olefinic.

## 4. MODELING APPROACH

The simulations were carried out using CONVERGE CFD software (version 3.0) [46]. Figure 2 presents the CFD domain and mesh used for the CFD simulation. In order to reduce the computational cost, only 1/7$^{th}$ of the single-cylinder geometry was simulated. The entire domain was uniformly meshed with a base grid of 1.4mm. In addition, the grid was refined in regions of interest using adaptive mesh refinement (AMR) and fixed embedding. Two levels of fixed embedding was adopted near the nozzle region and one level of fixed embedding near the walls of the cylinder. In addition, two levels of AMR were employed based on sub-grid criteria of 1m/s and 2.5K for velocity and temperature. Lastly, one level of AMR was employed to track the evaporation of each fuel component.

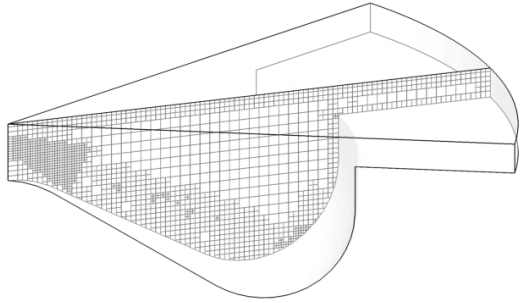

Figure 2: CFD domain simulated along with mesh resolution on a plane cut through center of injector.

In-cylinder turbulence was modeled using Reynolds-Averaged Navier Stokes (RANS)-based re-normalized group (RNG) $k$-$\varepsilon$ model [47] with wall functions. The wall heat transfer was accounted for by using the Han and Reitz model [48]. A second order-accurate spatial discretization scheme was used for the governing equations with a fully implicit first order-accurate time integration scheme. The liquid spray evolution processes were well characterized using state-of-the-art models available within CONVERGE. The injection was modeled using blob model [49] and the spray break-up processes were modeled using Kelvin-Helmholtz (KH) and Rayleigh-Taylor (RT) instability mechanisms [50]. The injector included angle (umbrella angle) and the cone angle of the parcels injected were 130º and 15º, respectively. Further, the interaction between liquid droplets and solid surfaces was accounted for by using the Wall Film model by O'Rourke and Amsden model [51]. Additional details about the CFD model and validation of the current spray model are discussed in detail in the authors' previous work [20].

In-cylinder combustion was modeled in each cell of the computational domain using the well-stirred reactor (WSR) combustion model. The finite-rate chemistry in each cell was modeled using a skeletal mechanism containing 241 species and 1396 reactions [20]. Further, the combustion model was coupled with a comprehensive soot model – Particle Mimic (PM) model based on Hybrid Method of Moments [52]. The physical and chemical processes accounted for in the soot model include – nucleation or particle inception through cyclopentafusedpyrene (A4R5), coagulation, surface growth through acetylene ($C_2H_2$), surface condensation through PAH species, and oxidation of soot by $O_2$ and OH.

# 5. RESULTS AND DISCUSSION
## 5.1 CFD MODEL VALIDATION

In previous work [20], the CFD model coupled with the above-mentioned skeletal mechanism was validated against experimental data in terms of in-cylinder pressure, apparent heat release rate (AHRR), CA10, CA50, and soot emissions for RD5-87 at four SOI timings: -21, -27, -36, and -45°aTDC. In the present study, the CFD model is further validated for two additional fuels – E30 and E100 in terms of various engine performance and emissions metrics. Figure 3 presents the comparison of in-cylinder average pressure and AHRR for E30 and E100 fuels at different SOI timings. The experimental pressure trace is an average of 300 consecutive cycles, while simulation pressure trace is based on the closed cycle simulation from IVC to EVO. It can be seen that the CFD model predicts the pressure and AHRR profiles reasonably well. However, the AHRR is over-predicted in simulations which is consistent with previous studies that adopted the Well Stirred Reactor (WSR) combustion model [20,53,54]. Further, the numerical model is compared with experimental data in terms of combustion phasing, soot, and $NO_x$ emissions. Figure 4 presents the comparison of CFD and experimental data in terms of CA10 (the crank angle at which 10% of fuel mass is consumed) and CA50 (the crank angle at which 50% of fuel mass is consumed) for E10, E30 and E100 fuels at four different SOI timings. The CFD model accurately predicts the CA10 and CA50 for all fuels within an uncertainty of one crank angle degree. Figure 5 (a&b) presents the comparison of soot mass and engine-out $NO_x$ emissions ($EINO_x$) emissions from experiments and simulations for E10, E30 and E100[2] at four different SOI timings. The non-monotonic trend in soot and $NO_x$ emissions with SOI timing for all fuels is well captured in the simulations. Soot emissions from CFD simulations capture the qualitative trend observed in experiments while quantitatively predicted within a factor of 2. On the other hand, the trends in NOx emissions for all fuels and with respect to SOI timing are preserved with marginal underprediction quantitatively. Overall, the validation study shows reasonable predictive capability of the CFD model and hence was adopted for further analysis.

---

[2] Soot emissions from E100 fuel were significantly low and within the uncertainty limit of the experimental setup, therefore not reported in experiments.

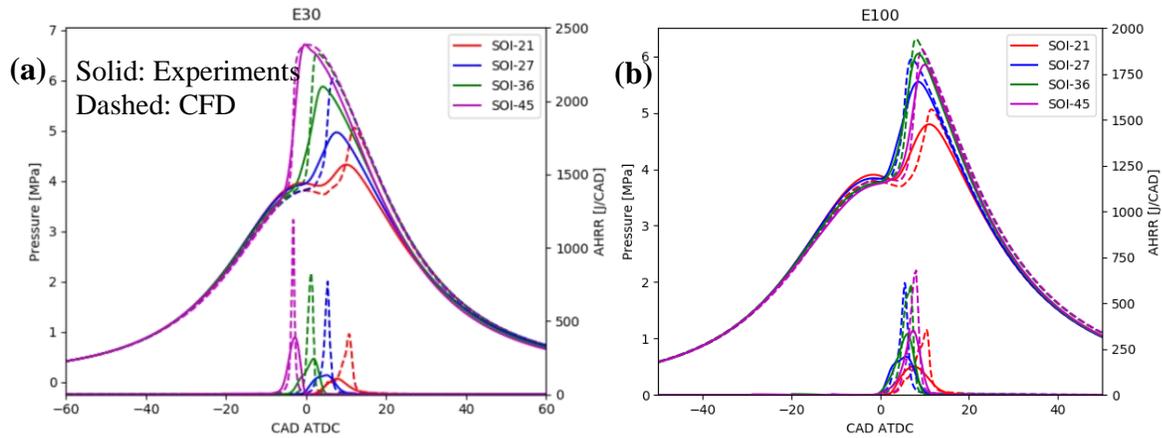

Figure 3: Comparison of CFD predictions against experimental data in terms of in-cylinder pressure and apparent heat release rates for (a) E30 and (b) E100 fuels.

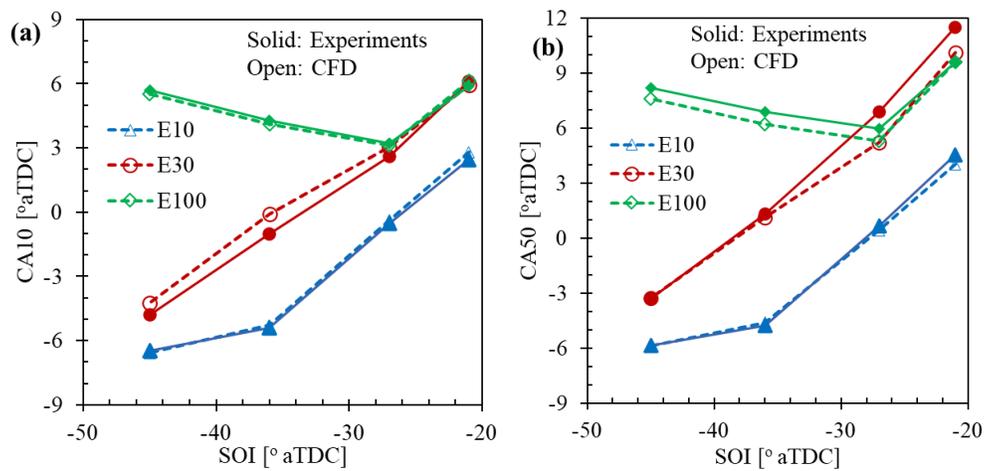

Figure 4: Comparison of CFD model predictions with experimental data in terms of (a) CA10 and (b) CA50.

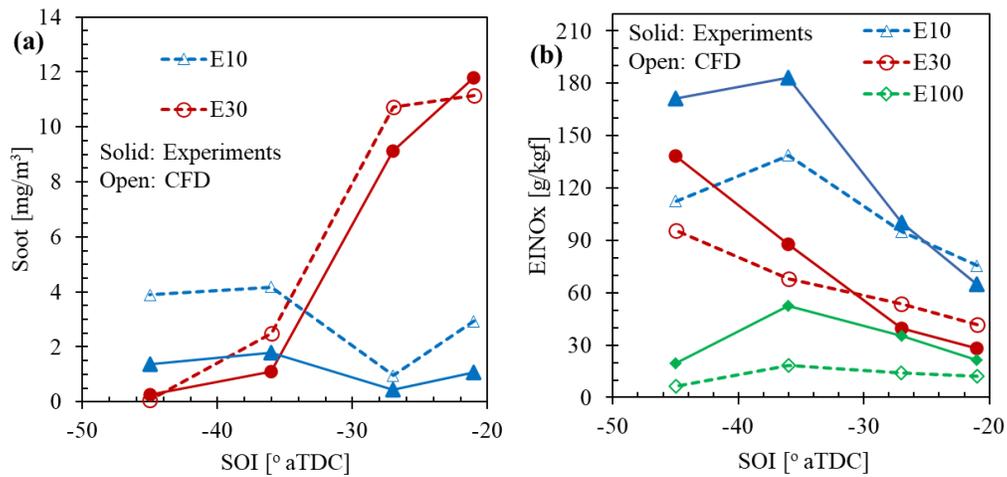

Figure 5: Comparison of results from experiments and simulations in terms of (a) soot mass and (b) EINO$_x$.

## 5.2 EFFECT OF FUEL COMPOSITION ON AUTOIGNITION

This section broadly discusses the effect of fuel composition (ethanol content) on autoignition propensity (represented by CA10). An attempt is made to isolate the effect of chemical reactivity from physical properties on autoignition and combustion phasing.

Figure 6 (a-b) presents the CA10 and ignition delay times[3] (IDT) for E10, E30, and E100 fuels at all SOI timings. In general, it is observed that E10 is more reactive than E30, followed by E100 at all SOI timings under consideration. Further, IDT increases with advancing SOI timing. In other words, the mixture becomes hard to ignite on advancing the SOI timing. This is a consequence of leaner mixtures at lower pressures and temperatures associated with advanced SOI timings, which make the in-cylinder conditions less favorable for autoignition. Table 3 presents the mean in-cylinder temperatures and pressures at SOI crank angle for all SOI timings. Lower temperatures lead to slower evaporation of liquid fuel resulting in poor fuel-air mixing process, thereby taking a longer time to establish conditions suitable for autoignition. As a consequence, a significant amount of fuel happens to interact and settle on the walls of the piston and liner with less amount of fuel available for combustion at advanced SOI timings. On the other hand, it is noticed that the rate of increase in IDT from -21 to -45°aTDC SOI timings is significantly higher for E100, while it is similar for E10 and E30 fuels. This is elaborated below in the context of isolating the effects of chemical reactivity and physical properties on autoignition propensity.

Table 3: In-cylinder mean pressure and temperature at crank angle corresponding to SOI timing.

| SOI [°aTDC] | Pressure (bar) | Temperature (K) |
|---|---|---|
| -21 | 26.4 | 976 |
| -27 | 20.0 | 930 |
| -36 | 13.1 | 844 |
| -45 | 8.60 | 777 |

---

[3] Ignition delay time (IDT) is estimated as the number of crank angles it takes to reach CA10 crank angle from start of injection (SOI) timing.

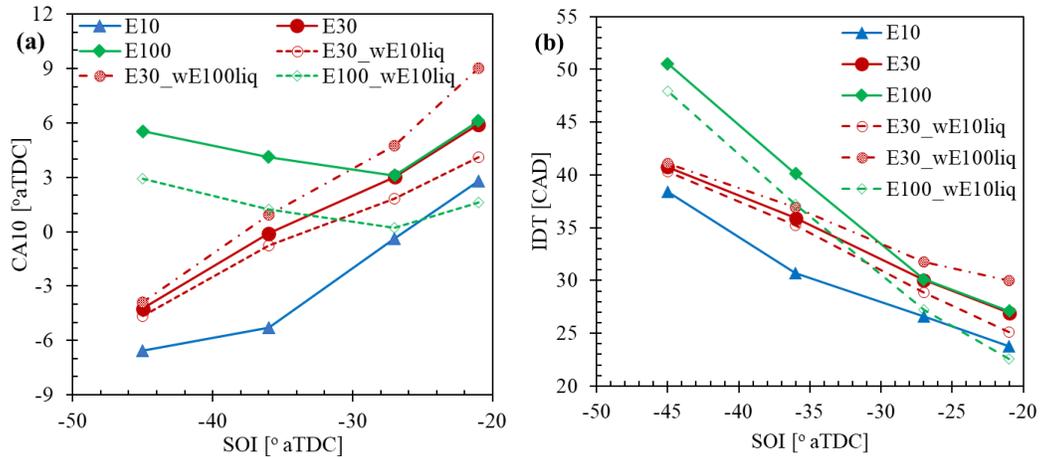

Figure 6: CA10 and ignition delay time (IDT) for E10, E30 and E100 fuels at SOI timings: -21/-27/-36/-45ºaTDC.

The IDT values follow the order E10<E30<E100 for all SOI timings. To better understand the reason for this behavior, an attempt is made to isolate the effect of physical properties (which have significant effect on the fuel evaporation process) from that of chemical reactivity. Thus, additional simulations were carried out with – E30_wE10liq, E30_wE100liq, and E100_wE10liq. To elaborate, E30_wE10liq case considers E30 fuel evaporating with physical properties[4] of E10 fuel while chemical reactivity is preserved to be that of E30. Comparing the results from E30 case with that of E30_wE10liq gives an understanding about the impact of E30's physical properties. On the other hand, comparing E30_wE10liq with E10 fuel case allows to understand the differences in chemical reactivity of E10, E30.

### 5.2.1 E10 vs E30

As discussed above, E10 is more reactive than E30 at all SOI timings. Comparing E30 with E30_wE10liq indicates that the autoignition is only marginally affected by the physical properties of E30. The change in CA10 with change in physical properties (E30 vs E30_wE10liq) is ≤ 1.5 CAD for all SOI timings. In other words, the differences in CA10 between E10 and E30 are mainly a consequence of differences in chemical reactivity. Figure 7 presents the ignition delay times for E10, E30, and E100 fuels at 25 bar pressure condition estimated from 0-D homogeneous constant volume model in CHEMKIN-Pro

---

[4] Physical properties include – density, viscosity, heat of vaporization, specific heat, surface tension, vapor pressure and conductivity.

[55]. The IDT in the 0-D simulations is defined as the time at which the temperature increases by 400K from the initial temperature. Further, the relative differences in 0D-IDTs between E10, E30, and E100 fuels remain similar at different pressures and hence are not presented. E10 and E30 share similar ignition delay times for temperatures >1000K (for all equivalence ratios relevant to the current conditions). However, for temperatures <1000K, E30 has higher IDT than that of E10 at all equivalence ratios. In addition, the IDTs for E30 are relatively more sensitive to change in temperature compared to E10 fuel. Further, the higher HoV of E30 fuel contributes to reducing in-cylinder temperatures when compared to E10. This is illustrated in Figure S1, which presents the density-weighted equivalence ratio and temperature in the cylinder for E10, E30, and E100 fuels at all SOI timings. These represent the equivalence ratio and temperature in the region of higher fuel vapor concentration in the cylinder [56]. The values are calculated based on Eqs. 1-2 at crank angle corresponding to CA5 of E10 fuel (since E10 is most reactive fuel).

$$\overline{\Phi} = \frac{\int \rho_F \Phi dV}{\int \rho_F dV} , \overline{T} = \frac{\int \rho_F T dV}{\int \rho_F dV} \qquad (1)$$

$$\rho_F = Y_F \times \rho_{total} \qquad (2)$$

where $\rho_F$ is local fuel vapor density excluding the liquid phase, $\rho_{total}$ is the total vapor density and $Y_F$ is mass fraction of fuel.

It is noticed that the density-weighted in-cylinder temperatures follow order E10>E30>E100. Thus, the lower reactivity of E30 at low temperatures when compared to E10 retards combustion for E30 at all SOI timings. Further, it was noticed that the differences in CA10 between E10 and E30 increase with advancing SOI timing. The differences increase from 3.1 at -21°aTDC to 5.2 CAD at -36°aTDC SOI timing. As previously discussed, advancing SOI timing results in cooler in-cylinder conditions, where E30 is significantly less reactive than E10 fuel (as shown in Figure 7a). However, in contrast, the difference in CA10 between E10 and E30 fuels at -45°aTDC SOI timing is slightly reduced to 2.8 CAD. This can potentially be a consequence of liquid fuel impinging on the lip of the piston which splits the amount of fuel injected between the piston bowl and squish region at -45°aTDC SOI timing [20]. The fuel in the squish region is not readily available for reaction resulting in poor conditions for autoignition. In addition, the low in-cylinder temperatures result in a significant amount of fuel to settle on the piston wall as wall film. Wall films further reduce the amount of fuel available for reaction. Overall, both of these phenomena make the autoignition process significantly complex, and thus difficult to justify the slightly reduced difference in CA10 between E10 and E30 fuels at -45°aTDC SOI timing.

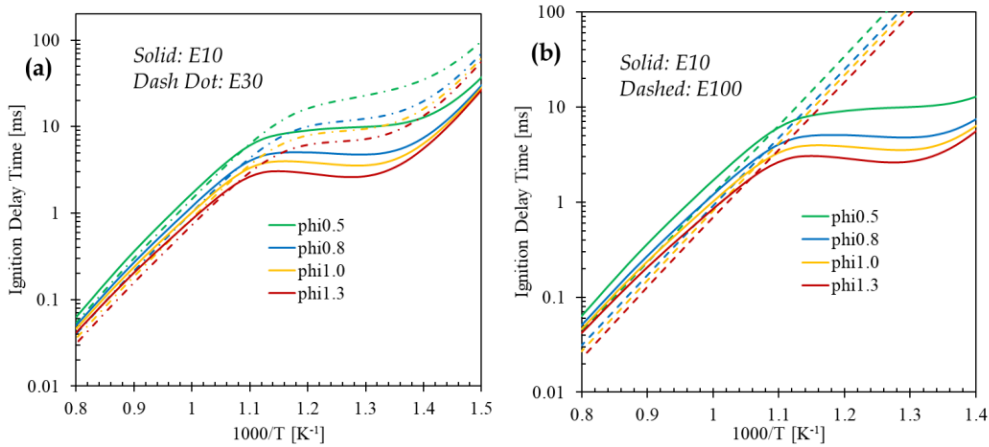

Figure 7: Comparison of ignition delay time between (a) E10 and E30, (b) E10 and E100 fuels.

### 5.2.2 E10 vs E100

E10 is more reactive (earlier CA10) when compared to E100 at all SOI timings. The differences in CA10 between E10 and E100 are due to the combined effect of physical properties and chemical reactivity, with the latter having a more significant effect at advanced SOI timings (-36 and -45°aTDC). The differences in CA10 between E10 and E100 are 3.2, 3.5, 9.4, and 12.1 CAD at -21, -27, -36, and -45°aTDC SOI timings, respectively. Comparing E10 with E100_wE10liq case allows to isolate the effect of chemical reactivity on autoignition. Swapping the physical properties of E100 with that of E10 advances CA10 of E100 fuel at all SOI timings. In other words, the physical properties of E100 make the fuel less reactive while swapping them with that of E10 allows to increase the reactivity. This is primarily associated with higher heat of vaporization (HoV) of E100 (795 kJ/kg @450K) than E10 (327 kJ/kg @450K) which lowers the in-cylinder temperature making the conditions less favorable for autoignition.

At -21 and -27°aTDC SOI timings, E10 and E100_wE10liq have similar CA10 (or IDT). This trend is consistent with results from 0D-IDT simulations (shown in Figure 7b). For temperatures > 900K, E10 and E100 share similar reactivity or similar IDTs (at all equivalence ratios). These conditions are relevant to -21 and -27°aTDC SOI timings and hence we observe similar CA10 for E10 and E100_wE10liq cases. Similarly, CA10 for -36 and -45°aTDC SOI conditions are advanced when physical properties of E100 are swapped with E10. However, the autoignition timing (CA10) for E100_wE10liq is still significantly delayed than E10 at these advanced SOI timings. These advanced SOI timings are associated with low in-cylinder temperatures (~850K) at which E100 has high resistance to autoignition (as shown in Figure 7b). The 0D-IDTs for E100 are substantially higher than E10 for temperatures < 900K. Moreover, the 0D-IDTs for E100 increases exponentially with temperature at all equivalence ratios, indicating very high sensitivity to temperature. As a result, the CA10 of E100 fuel for -36 and -45°aTDC SOI cases occur at significantly retarded crank angles, which happen to be closer to crank angles corresponding to CA10 of -21 and -27°aTDC SOI timings. Overall, the chemical reactivity of the fuel has a much stronger effect on autoignition timing when compared to physical properties.

The differences in CA10 between E30 and E100 at different SOI timings are qualitatively consistent with differences in CA10 between E10 and E100. Hence, the differences in CA10 between E30 and E100 are not discussed for the sake of brevity.

**5.3 EFFECT OF FUEL COMPOSITION ON SOOT EMISSIONS**

This section broadly discusses soot emissions for E10, E30 and E100 fuels. Figure 8 presents the soot mass for E10, E30 and E100 fuels at four SOI timings. It is noticed that E30 produces significantly more soot than E10 for retarded SOI timings (-21 and -27°aTDC). In contrast, E30 produces less soot than E10 for -36 and -45°aTDC SOI conditions. On the other hand, the trend in soot mass produced from E10 follows an order: -36>-45>-21>-27°aTDC SOI, while for E30 soot mass monotonically decreases with advancing SOI timing. Soot formation under the low-load GCI conditions considered here is primarily a consequence of fuel film accumulated on piston walls (resulting in pool fires). Low in-cylinder temperatures due to advanced SOI timings and high HoV of gasoline/ethanol blends makes it unfavorable for evaporation of liquid fuel, resulting in impingement of liquid fuel on piston walls. As a consequence, we notice significantly high amounts of fuel film accumulation on piston walls. The film accumulated evaporates as combustion progresses, however oxygen deficient conditions result in the formation of soot precursors (acetylene, ethylene, benzene, and other larger polycyclic aromatic hydrocarbons) and thus, soot. Overall, film accumulation and the rate of evaporation of film play a significant role in the formation of soot under these GCI conditions. This phenomenon has been discussed in detail in our previous publication in the context of gasoline/ethanol/n-butanol blends [20]. Further, the details about soot mass produced from E10 and the trend observed concerning SOI timings were also elaborated in our previous publication [20]. Hence, the current section focuses on discussing the sooting tendency of E30 and E100 fuels and their differences relative to E10 fuel.

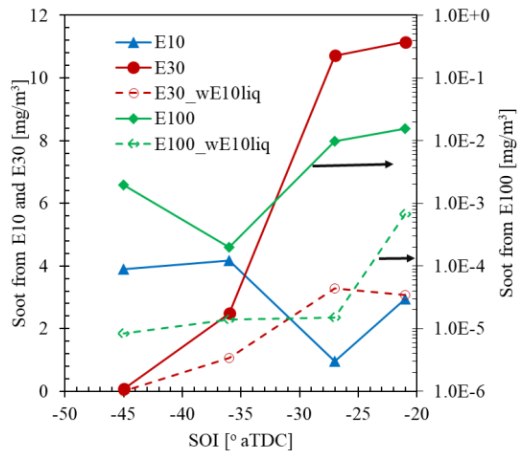

Figure 8: Soot mass for E10, E30, and E100 fuels at four SOI timings.

### 5.3.1 E10 vs E30

Figure 9 presents the ratio of fuel film mass on piston wall to fuel mass injected from E10, E30, and E30_wE10liq. The amount of film formed follows the order E30>E10≥E30_wE10liq for all SOI timings. Comparing soot mass from E30 with that from E30_wE10liq allows for characterization of the impact of physical properties of E30 on soot emissions. E30_wE10liq produces lower soot emissions when compared to E30 at all SOI timings. This indicates that physical properties of E30 increase the sooting tendency of E30 fuel. This is a consequence of significantly higher film mass accumulation for E30 fuel when compared to E30_wE10liq. For instance, -21°aTDC SOI condition has ~9.5% and 5.5% of fuel mass accumulated as film on piston walls at TDC for E30 and E30_wE10liq cases, respectively. Figure 10 shows the contours of film mass and soot mass at various crank angles for SOI-21 corresponding to E30 and E30_wE10liq cases. It is noticed that soot formation starts close to CA50 (10.1 and 7.2º for E30 and E30_wE10liq, respectively) for both cases and happens close to the piston wall. As previously mentioned, the film starts to evaporate and results in the formation of soot. Lower film mass for E30_wE10liq results in lower soot emissions. To elaborate, the higher film mass accumulated for E30 when compared to E30_wE10liq is mainly a consequence of higher HoV of E30 compared to E10. The HoV of E10 and E30 are 327 and 508 kJ/kg at 450K. In other words, the higher HoV of E30 requires higher amount of heat to evaporate resulting in cooler in-cylinder conditions. As a consequence, the evaporation of remaining liquid fuel slows down leading to impingement of liquid fuel on the walls of the cylinder. Figure S4 presents the ratio of film mass and injected fuel mass from E30 fuel when adopted with different combinations of E10 and E30 physical properties. wE10Cp represents the case where only the specific heat of E30 fuel is swapped with that of E10 while other physical properties remain same as those of E30. wE10 indicates the case same as E30_wE10liq. Comparing the effect of different physical properties, it can be concluded that high HoV of E30 leads to significantly high film mass accumulation at the piston wall thus, resulting in higher soot emissions. It is important to note, however, that while E30 produces lower soot emissions than E30_wE10liq at all SOI timings, the effect of physical properties diminishes with advance in SOI timing.

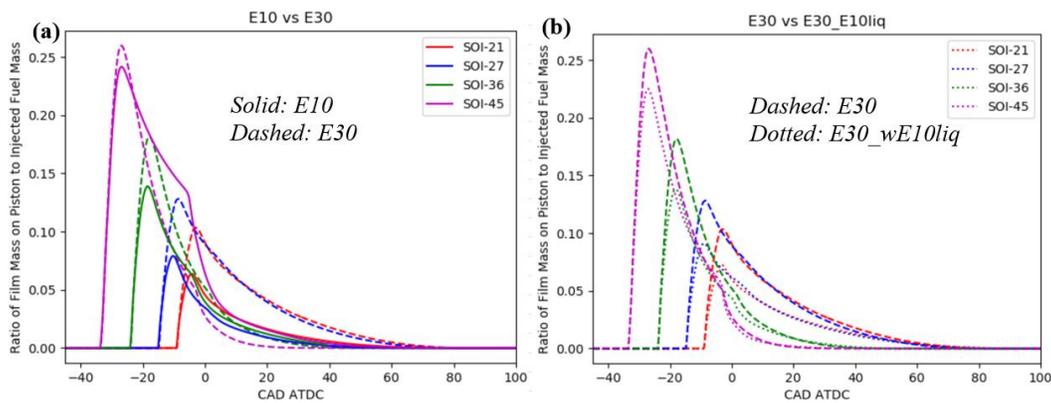

Figure 9: Ratio of film mass on piston to injected fuel mass plotted versus CAD for (a) E10 and E30, (b) E30 and E30_wE10liq cases.

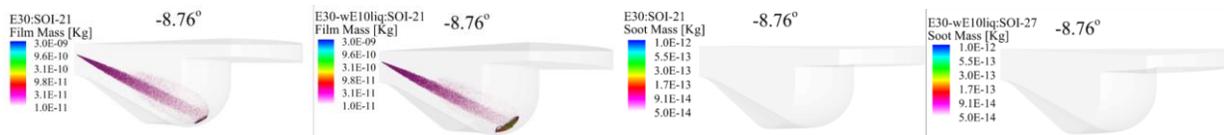

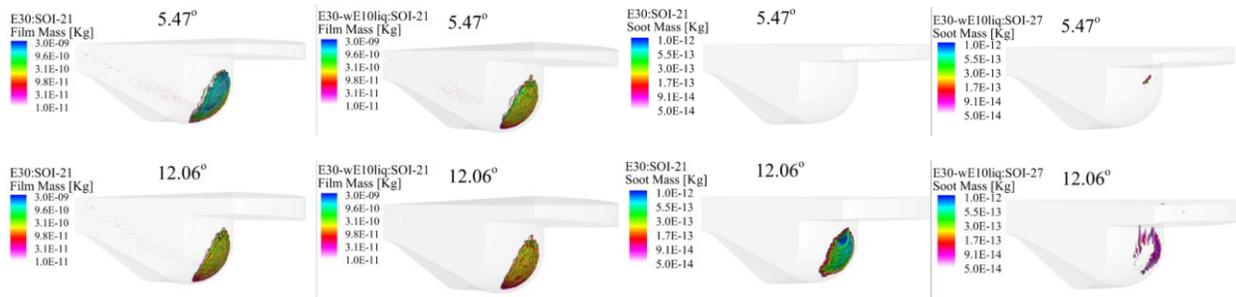

Figure 10: Contours showing film mass and soot mass at different crank angles for E30 and E30_wE10liq cases corresponding to -21°aTDC SOI timing.

Comparing E10 with E30_wE10liq allows for characterizing the chemical effect of E10 and E30 on soot emissions. As indicated in Fig. 8, E30_wE10liq exhibits marginally higher sooting tendency than E10 at retarded (-21 and -27°aTDC) SOI timings while it becomes less sooting at advanced SOI timings (-36 and -45°aTDC). From fundamental understanding, E30 is less sooting due to the lower aromatic content (toluene) in the fuel and higher oxygen concentration in the fuel when compared to E10 fuel. In TPRF-alcohol fuels, the formation of soot precursors has a strong correlation with toluene content of the fuel [31, 35]. Higher toluene content in the fuel blend results in higher concentration of soot precursors and thus soot. Further, the higher oxygen content in the fuel can potentially aid in the oxidation of soot precursors and thus less soot formed. Thus, lower toluene and higher oxygen content aid in reducing the sooting tendency of E30 fuel when compared to E10, assuming these fuels have the same physical properties. Figure S5 shows the mass fraction of soot precursors - acetylene, pyrene and cyclopentafusedpyrene from E10 and E30_wE10liq cases at -45°aTDC SOI timing. It is evident that the amount of soot precursors produced with E30_wE10liq is significantly lower than that produced with E10 fuel. Hence, at advanced SOI timings where E10 and E30_wE10liq have similar film mass accumulation at piston wall, E30_wE10liq produces lower soot. This is further confirmed by Figure 11 which presents contours of OH mass fraction and soot formation in E10 and E30_wE10liq cases at -45°aTDC SOI timing. The mass fraction of OH is higher for E30_wE10liq while soot formation (defined as sum of contributions from nucleation, coagulation, condensation, and segregation) is higher for E10. Consequently, soot precursors get oxidized resulting in lower soot formation in E30 blend. Moreover, the higher concentration of OH also promotes the oxidation of soot, thereby resulting in significantly lower net soot emissions from E30 at advanced SOI timings.

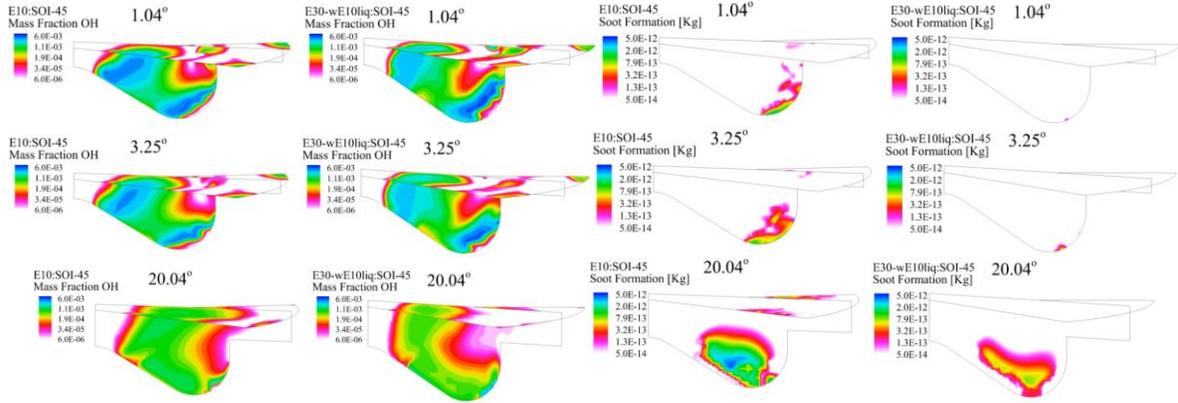

Figure 11: Contours of mass fraction of OH mass fraction and soot mass for E10 and E30_wE10liq cases at -45°aTDC SOI timing.

On the other hand, the anomaly observed at retarded SOI timings is mainly due to the lower reactivity of E30 when compared to E10. CA10 and CA50 values for E30_wE10liq are significantly lower compared to those for E10. For instance, at -21°aTDC SOI timing, CA10, CA50 for E10 and E30_wE10liq are 2.8, 4.03 and 4.2 and 7.3°aTDC, respectively. The retarded combustion phasing limits the oxidation of soot formed due to lower residence time. As a consequence, although E30_wE10liq has higher OH concentration than E10, we notice higher soot concentration for retarded SOI timings. On the other hand, although E30_wE10liq has retarded combustion phasing than E10 at advanced SOI timings, significantly higher OH concentration and higher residence time available allows for efficient oxidation of soot precursors and soot formed. Thus, E30_wE10liq has lower soot emissions than E10 at advanced SOI timings.

Overall, the lower toluene content and higher ethanol concentration (that can yield more OH radicals) make E30 chemically less prone to produce soot than E10. However, the physical properties (mainly HoV) make E30 significantly more sooting fuel than E10 at retarded SOI timings (-21 and -27°aTDC). On the other hand, the chemical effect dominates for E30 making it less sooting than E10 at advanced SOI timings.

### 5.3.2 E10 vs E100

Sooting propensity of E100 fuel is significantly lower than those of E10 and E30 at all SOI timings. Figure 12 presents the ratio of film mass on piston walls to injected fuel mass for E10 and E100 fuels at all SOI conditions considered. It is noticed that the fuel mass accumulated on the piston walls is significantly higher for E100 fuel when compared to E10 (and also E30 although not shown). However, the higher fuel film only affects the combustion phasing by reducing the availability of fuel for combustion. The lower carbon content and higher oxygen content in E100 significantly reduce the sooting propensity of E100 fuel. This can be confirmed by comparing soot emissions from E100_wE10liq and E10 cases. It is evident from Figure 8 that soot from E100_wE10liq is lower than E10 by a factor of ~0.5E5 at all SOI timings, thus clearly indicating the low sooting propensity of E100 fuel at all SOI timings. Considering the very low soot emissions from E100, further investigation on comparing soot from E100 and E10 soot is not considered.

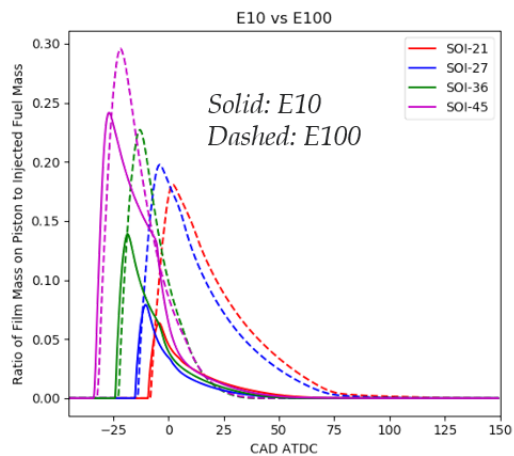

Figure 12: Ratio of film mass on piston to injected fuel mass from E10 and E100 fuels at four SOI timings.

## 5.4 EFFECT OF SOI TIMING AND FUEL COMPOSITION ON $NO_x$ EMISSIONS

This section discusses the effect of fuel composition on $NO_x$[5] emissions at different SOI timings. Figure 13 presents the NOx emission index versus SOI for E10, E30, and E100 fuels. It is noticed that $NO_x$ emissions follow the order E10>E30>E100 at all SOI timings. Formation of $NO_x$ is highly dependent on the maximum temperature of the burning gases, oxygen content, and residence time available for the reactions to occur at these extreme conditions [57,58,59]. This section is split into two sub-sections to further characterize the effect of SOI timing and fuel composition.

---

[5] $NO_x$ is defined as sum of concentrations of NO and $NO_2$.

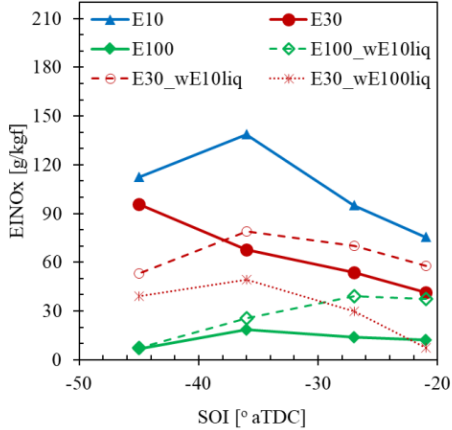

Figure 13: EINO$_x$ for E10, E30, and E100 fuels at -21/-27/-36/-45°aTDC SOI timings.

### 5.4.1 Effect of SOI timing on NO$_x$ emissions

As indicated earlier, EINO$_x$ increases with advancement in SOI timing from -21 to -36°aTDC. Retarded injection timing allows fuel to be introduced into the combustion chamber closer to the top dead center (TDC), when the in-cylinder temperatures and pressure are relatively high. This eventually leads to rapid fuel evaporation, cooling of the air-fuel mixture and reduction in peak combustion temperatures. On the other hand, the mixing time available is relatively lower at retarded SOI timings. Consequently, the fuel-air mixture is more stratified at retarded SOI timings (-21 and -27°aTDC), which causes autoignition to occur at multiple locations. This increases in-cylinder temperatures steadily allowing to burn the rest of the mixture. On the other hand, advanced SOI timings allow more time for fuel-air mixing leading to homogeneous mixtures with lower equivalence ratios. With lower mixture stratification, the ignition tends to occur in a bulk (volumetric) fashion resulting in steep rise in pressure and temperatures. Further details about the mixture stratification and autoignition phenomena for E10 fuel are discussed in our previous publication [20]. Moreover, advanced SOI timings happen to have autoignition occurring well before TDC which allows the in-cylinder temperature to further increase, resulting in higher in-cylinder temperatures than retarded SOI timings where ignition only occurs after TDC during expansion stroke. Thus, the increase in NO$_x$ emissions at advanced SOI timing is a consequence of the mixture being more homogeneous and the autoignition occurring well before TDC that results in higher in-cylinder temperatures. However, NO$_x$ emissions decrease marginally with advancing SOI timing from -36 to -45°aTDC for E10 and E100 fuels.

At -45°aTDC SOI timing, the spray impinges at the lip of the piston splitting the fuel between piston and squish regions. Fuel in the squish region is exposed to higher heat transfer rates making it not readily available for combustion, resulting in slower rise in in-cylinder temperatures. Figure 14 presents the average in-cylinder temperatures for E10, E30, and E100 fuels at all SOI timings. It is evident that the temperatures for -45°aTDC SOI are lower than those for -36°aTDC SOI for both E10 and E100 fuels (Figure 14a). In addition, for E100 fuel the decrease in temperature is also due to very low reactivity of E100 at this advanced injection timing which significantly retards autoignition (as discussed in section 5.1) beyond TDC while for E10 and E30, autoignition for -45°aTDC SOI timing occurs well before TDC. As a consequence, although autoignition allows to increase temperature by burning fuel-air mixture, the in-cylinder mixture temperature tends to decrease faster due to the expansion stroke. This results in reduced potential for thermal $NO_x$ production. Thus, we notice a decrease in $NO_x$ emissions when SOI timing is advanced beyond -36°aTDC. In contrast, $NO_x$ emissions increase for E30 fuel when SOI is advanced from -36 to -45°aTDC. This is associated with higher in-cylinder temperatures for -45°aTDC SOI condition when compared to -36°aTDC SOI condition (as shown in Figure 14b). The potential reason for higher temperatures is associated with faster evaporation of fuel film accumulated on the piston for E30 when compared to E10 (as shown in Figure 9a), making more fuel readily available for combustion when compared to E10. This is also evident from reduced difference in IDT when advancing SOI timing from -36 and -45°aTDC for E30 when compared to E10 (as shown in Figure 6b). In other words, the CA10 advances by 1.3° and 4.1° for E10 and E30, respectively when SOI timing is advanced from -36 to -45°aTDC. Earlier combustion implies the potential for higher in-cylinder temperatures, which results in higher $NO_x$ emissions. Thus, higher in-cylinder temperatures for E30 fuel when advancing SOI from -36 to -45°aTDC results in an increase in $NO_x$ emissions.

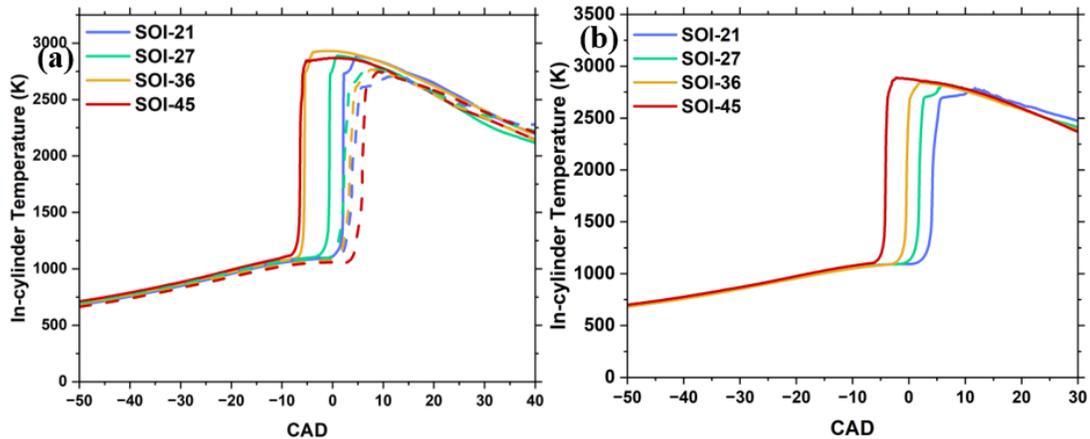

Figure 14: In-cylinder average temperatures for (a) E10 (solid) and E100 (dashed) (b) E30 fuel at -21/-27/-36/-45°aTDC SOI timings.

### 5.4.2 Effect of fuel composition on NO$_x$ emissions

NO$_x$ emissions follow the trend E10>E30>E100 for all SOI timings. The trend is consistent with the trend observed in autoignition propensity of E10, E30, and E100 fuels. As previously discussed, retarded combustion phasing is associated with lower in-cylinder temperatures and consequently lower thermal NO$_x$ emissions. Changing the physical properties of E30 and E100 with that of E10 increases NO$_x$ emissions marginally (as shown by E30_wE10liq and E100_wE10liq cases in Figure 13). This increase in NO$_x$ with changing physical properties to E10 is mainly due to the lower HoV of E10 fuel which allows to keep the in-cylinder conditions hotter than E30 and E100. However, NO$_x$ emissions from E30_wE10liq and E100_wE10liq are still lower than those from E10, which indicates there is an impact of differences in chemical reactivity of fuels on NO$_x$. Overall, the results remain consistent with trends observed in autoignition propensity of E30_wE10liq and E100_wE10liq. Thus, it can be concluded that the trend observed in NO$_x$ emissions from E10, E30, and E100 fuels is due to a combination of physical properties (mainly HoV) and reactivity order of fuel – E10>E30>E100 at conditions investigated in this study.

### 6. CONCLUSIONS

The present computational study investigated the effects of SOI (start of injection) timing, and the physical and chemical properties of three gasoline-ethanol blends with varying ethanol content, on autoignition propensity, combustion phasing, soot and NO$_x$ emissions under low-load GCI engine conditions. The fuels considered are E10, E30, and E100 containing 10, 30 and 100% of ethanol by volume in the fuel mixture, respectively. A four-component surrogate TPRF + Ethanol was considered to capture the chemical and physical behavior of the fuels considered. The experiments were conducted in a Heavy-Duty Caterpillar (CAT) 3401 single-cylinder engine upgraded for common rail injection. The same engine was modeled using advanced computational fluid dynamic (CFD) modeling tools. The CFD model was coupled with a full chemistry solver that used a detailed chemical kinetic mechanism, containing 241 species. The coupled model was validated against experimental data in terms of in-cylinder pressure, heat release rate, $NO_x$, and soot emissions for E10, E30, and E100 fuels at four SOI timings (-21/-27/-36/-45°aTDC). The model captured the experimental data both qualitatively and quantitatively. Important observations are as follows.

1. The ignition delay time (IDT) of all fuels increased with advancing SOI timing, while the increase in IDT for E100 fuel was significantly higher for advanced SOI timings (-36 and -45°aTDC) when compared to E10 and E30. The high HoV of E100 reduced the in-cylinder temperatures significantly and the very low reactivity of E100 at low temperatures increased the IDT for E100 fuel at -36 and -45°aTDC SOI timings.

2. The autoignition propensity of fuels follows the order E10>E30>E100 with E10 having the earliest CA10 at all SOI conditions. It was observed that the chemical reactivity of fuels had the dominant effect on autoignition propensity at all SOI timings. In contrast, the fuel physical properties (density, specific heat, viscosity, heat of vaporization, conductivity, surface tension and vapor pressure) had no significant effect on CA10 at advanced SOI timings (-36 and -45°aTDC), and a moderate effect at retarded SOI timings (-21 and -27°aTDC).

3. The sooting propensity of E100 was lower by a factor of more than 100 compared to that of E10 and E30 at all SOI timings. Thus, E100 can be considered soot free under the operating conditions considered in this study.

4. E30 fuel, having higher ethanol content produced lesser soot emissions than E10 at earlier SOI timings. However, E30 produced significantly higher soot than E10 at retarded SOI timings (-21 and -

27ºaTDC). While E30 is a chemically less sooting fuel, its physical properties cause the fuel to be more sooting than E10 at retarded SOI timings. This is a consequence of retarded autoignition (that reduces residence time available to oxidize soot) and high fuel-film accumulation on walls.

5. $NO_x$ emissions followed the order E10>E30>E100 at all SOI timings. This was consistent with the autoignition propensity of the fuels. Retarded combustion phasing resulted in lower in-cylinder temperatures, leading to less favorable conditions for $NO_x$ formation.

6. Overall E100 has comparable reactivity as E30 at late injection timings (-21 and -27ºaTDC) with 100- and 5-times lower soot and $NO_x$ emissions, respectively, than E30.

# 7. ACKNOWLEDGEMENTS

The submitted manuscript has been created by UChicago Argonne, LLC, Operator of Argonne National Laboratory (Argonne). Argonne, a U.S. Department of Energy Office of Science laboratory, is operated under Contract No. DEAC02-06CH11357. The U.S. Government retains for itself, and others acting on its behalf, a paid-up nonexclusive, irrevocable worldwide license in said article to reproduce, prepare derivative works, distribute copies to the public, and perform publicly and display publicly, by or on behalf of the Government. This research was partially funded by DOE's Office of Vehicle Technologies and Office of Energy Efficiency and Renewable Energy under Contract No. DE-AC02-06CH11357. The authors wish to thank Gurpreet Singh, Kevin Stork, and Michael Weismiller, program managers at DOE, for their support. This research was conducted as part of the Co-Optimization of Fuels & Engines (Co-Optima) project sponsored by the U.S. Department of Energy (DOE) Office of Energy Efficiency and Renewable Energy (EERE), Bioenergy Technologies and Vehicle Technologies Offices. The authors would like to acknowledge the Laboratory Computing Resource Center (LCRC) at Argonne National Laboratory and computing facilities at the University of Illinois at Chicago (UIC), for computational resources used in this research. The authors wish to thank Russ Fitzgerald and many others at Caterpillar Inc. for support in upgrading Argonne's Caterpillar 3401 single-cylinder research engine to common rail operation, providing Caterpillar Next Gen injectors, and supplying information on their function to the CFD modeling in this work.